\documentclass[]{amsart}

\def\ba{\mathbb{A}}
\def\bi{\mathbb{I}}
\def\bj{\mathbb{J}}
\def\bk{\mathbb{K}}
\def\bM{\mathbb{M}}
\def\bn{\mathbb{N}}

\def\bc{\mathbb{C}}
\def\bl{\mathbb{L}}
\def\bs{\mathbb{S}}
\def\bff{\mathbb{F}}
\def\h{\mathcal{H}}
\def\k{\mathcal{K}}
\def\l{\mathcal{L}}
\def\g{\mathcal{G}}
\DeclareMathOperator{\B}{\rm B}

\newcommand{\zn}{Z^{\bn}}
\newcommand{\om}{\mathcal{O}}
\newcommand{\com}[1]{#1^{\prime}}
\newcommand{\ort}[1]{#1^{\perp}}
\newcommand{\norm}[1]{\Vert#1\Vert}
\newcommand{\cbnorm}[1]{\Vert#1\Vert_{{\rm cb}}}
\newcommand{\inner}[2]{\langle #1,#2\rangle}

\newcommand{\suma}{\sum_{\alpha\in\mathbb{A}}}
\newcommand{\sumi}{\sum_{i\in\mathbb{I}}}

\newcommand{\bb}[2]{{\rm B}(#1,#2)}
\newcommand{\bh}{{\rm B}(\mathcal{H})}
\newcommand{\kl}{{\rm K}_{\ell}}

\newcommand{\al}{{\rm A}_{\ell}}
\newcommand{\ar}{{\rm A}_{r}}
\newcommand{\pl}{\phi_{\ell}}
\newcommand{\bkh}{{\rm B}(\mathcal{K},\mathcal{H})}
\newcommand{\kkh}{{\rm K}(\mathcal{K},\mathcal{H})}
\newcommand{\bm}[3]{{\rm B}_{#1}(#2,#3)}
\newcommand{\kh}{{\rm K}(\mathcal{H})}
\newcommand{\acbbb}[2]{{\rm CB}_A(#1,#2)_B}
\newcommand{\cbb}[2]{{\rm CB}(#1,#2)}
\newcommand{\matn}[1]{{\rm M}_{n}(#1)}

\newcommand{\matj}[1]{{\rm M}_{\mathbb{J}}(#1)}
\newcommand{\matk}[1]{{\rm M}_{\mathbb{K}}(#1)}

\newcommand{\rowj}[1]{{\rm R}_{\mathbb{J}}(#1)}
\newcommand{\colj}[1]{{\rm C}_{\mathbb{J}}(#1)}
\newcommand{\mati}[1]{{\rm M}_{\mathbb{I}}(#1)}
\newcommand{\rowi}[1]{{\rm R}_{\mathbb{I}}(#1)}
\newcommand{\coli}[1]{{\rm C}_{\mathbb{I}}(#1)}
\newcommand{\aob}{{_A{\rm OM}_B}}
\newcommand{\mnon}{{_M{\rm NOM}_N}}
\newcommand{\ha}{{_A{\rm HM}}}
\newcommand{\hb}{{_B{\rm HM}}}
\newcommand{\nhm}{{_M{\rm NHM}}}
\newcommand{\nhn}{{_N{\rm NHM}}}
\newtheorem{theorem}{Theorem}[section]
\newtheorem{lemma}[theorem]{Lemma}
\newtheorem{corollary}[theorem]{Corollary}
\newtheorem{proposition}[theorem]{Proposition}
\theoremstyle{remark}
\newtheorem{remark}[theorem]{Remark}
\theoremstyle{definition}
\newtheorem{definition}[theorem]{Definition}
\newtheorem{example}[theorem]{Example}
\numberwithin{equation}{section}

\begin{document}

\title[]{Injective cogenerators among operator bimodules} 
\author{Bojan Magajna} 
\address{Department of Mathematics\\ University of Ljubljana\\
Jadranska 19\\ Ljubljana 1000\\ Slovenia}
\email{Bojan.Magajna@fmf.uni-lj.si}



\begin{abstract} Given C$^*$-algebras $A$ and $B$ acting cyclically on Hilbert spaces $\h$ 
and $\k$, respectively, we characterize completely isometric $A,B$-bimodule maps from 
$\bkh$ into  operator $A,B$-bimodules. We determine  cogenerators in some classes of 
operator bimodules. For an injective cogenerator $X$ in a suitable category of operator 
$A,B$-bimodules we show: if $A$, regarded as a C$^*$-subalgebra of  $\al(X)$ (adjointable 
left multipliers on $X$), is equal to its relative double commutant in $\al(X)$, then $A$ 
must be a W$^*$-algebra.
\end{abstract}

\maketitle

\section{Introduction}

An operator space $Z$ is called {\em injective} provided that for each inclusion of operator
spaces $X\subseteq Y$ every completely bounded map $\phi$ from $X$ to $Z$ extends to a completely
bounded map $\tilde{\phi}$ from $Y$ to $Z$ with the same completely bounded norm. If, 
in addition, $X,Y$ and $Z$ are operator bimodules over a pair of C$^*$-algebras $A$ 
and $B$ and $\phi$ is a bimodule map, then $\tilde{\phi}$ can be achieved to be an $A,B$-bimodule
map too, and $Z$ is then an {\em injective operator $A,B$-bimodule}. This property, shared
by $Z=\bh$ (by the well known extension theorem \cite{P}, \cite{W}), enables one to treat completely
bounded bimodule maps 
into $Z$, to a certain extent, as linear functionals. But there is another, equally
important, property of linear functionals that is not contained in the concept of 
injectivity, namely linear functionals separate points of the space. 
An operator $A,B$-bimodule $Z$ is called a {\em cogenerator} if for every operator 
$A,B$-bimodule $X$, every $x\in X$ and $\varepsilon>0$ there exists an $A,B$-bimodule complete contraction $\phi:X\to Z$
such that $\norm{\phi(x)}>\norm{x}-\varepsilon$. If (in addition) $\phi$
can be found such that $\norm{\phi(x)}=\norm{x}$, then $Z$ is called a {\em strict
cogenerator}. Not all injective operator bimodules are 
cogenerators: for example, it is not hard to show (but will not be needed here) that 
a continuous von Neumann algebra $A$ is not a cogenerator as an operator $A$-bimodule. 

Before describing  the content of the paper, 
let us mention briefly a wider context for motivation. In pure
algebra duality for modules has been 
intensively studied at least since the appearance of
paper  \cite{Mo} by Morita (see also  \cite[Chapter 19]{L} for a more modern account). 
While the related notion of Morita equivalence 
has been vigorously studied  also outside pure algebra, in particular 
in  operator algebra theory (see \cite{R},
\cite{BMP}, \cite{B2} and the references there), the operator module duality itself has been 
considered so far only for the range bimodules of
special nature (\cite{N}, \cite{Po}, \cite{M}). 
(We note at this point that, unlike in pure algebra, there is no natural reduction of 
operator $A,B$-bimodules 
to, say, left modules over $A\otimes B^{\rm op}$, since
there is in general no operator algebra norm structure known on $A\otimes B^{\rm op}$ 
turning all operator $A,B$-bimodules  into operator left 
$A\otimes B^{\rm op}$-modules.) However, the recent theory of multipliers of operator
spaces (\cite{B4}, \cite{BP}, \cite{BEZ}, \cite{We}) enables us a natural definition of 
duality.  Namely, according 
to Blecher \cite[5.4]{B4}
for each operator space $X$ 
there exist C$^*$-algebras $\al(X)$ and $\ar(X)$ 
such that $X$ is an operator $\al(X),\ar(X)$-bimodule and each operator $A,B$-bimodule 
structure on $X$  is induced by some
$*$-homomorphisms $A\to \al(X)$ and $B\to\ar(X)$. This suggests  
the following definition. 

\smallskip
Let $\alpha:A\to \al(X)$ and $\beta:B\to \ar(X)$ be $*$-homomorphisms  and denote by $A^c$ 
and $B^c$ the commutants
of $\alpha(A)$ and $\beta(B)$ in $\al(X)$ and $\ar(X)$, respectively. Then the
{\em dual of an operator $A,B$-bimodule $Y$ with respect to $X$} is the operator 
$A^c,B^c$-bimodule $\acbbb{Y}{X}$ (completely bounded $A,B$-bimodule
maps from $Y$ to $X$),  where
$(c\phi)(y):=c\phi(y)$ and  $(\phi d)(y):=\phi(y)d$ ($y\in Y,\ c\in A^c,\ d\in B^c$).

\smallskip
Restricting in this definition $X$ to be an (injective) cogenerator will guarantee the 
existence of enough
`functionals'. Then the bidual is an operator bimodule over C$^*$-algebras
$A^{cc}$ and $B^{cc}$, which contain $A$ and $B$, respectively. We shall say that the pair
$(A,B)$ {\em admits a duality} if for some injective cogenerator $X$ we have  $A^{cc}=A$
in $\al(X)$ and $B^{cc}=B$ in $A_r(X)$. (Here we have identified $A$ with $\alpha(A)$
and $B$ with $\beta(B)$.) In Section 5 we prove that $(A,B)$-admits a duality if and only if
$A$ and $B$ are von Neumann algebras. For this we shall need some preparation. 

Let $\h$ and $\k$ be Hilbert modules over C$^*$-algebras $A$ and $B$ (respectively), $X$
an operator $A,B$-bimodule and $\phi:\bkh\to X$ a completely contractive $A,B$-bimodule map.
In Section 3 we will prove that, if $\h$ and $\k$ are cyclic with unit cyclic vectors
$\xi\in\h$ and $\eta\in\k$ and $\norm{\phi(\xi\otimes\eta^*)}=1$,  
then $\phi$ is automatically completely isometric. In Section 4 we will show that if 
$X$ is injective, $\phi$ induces a completely contractive and completely positive
$A$-bimodule map $\pl:\bh\to\al(X)$,
which is completely isometric if $\phi$ is completely isometric and $\h$ is cyclic. 
This will imply that
$\al(X)$ is a cogenerator for operator $A$-bimodules if $X$ is a cogenerator for operator 
$A,B$-bimodules. Further, if  $\phi$ is injective, then
the kernel of $\pl$ does not contain positive (nonzero) operators. 

Now assume that
$A\subseteq\bh$,
$\pi$ is a $*$-homomorphism from $A$ into an injective C$^*$-algebra $C$ 
(say, $C=\al(X)$ when $X$ is injective) making $C$ an $A$-bimodule, and $\psi:\kh\to C$
is a bounded left $A$-module map, the kernel of which does not contain nonzero positive elements.
In this situation we will prove in Section 6 that, if $\pi(A)$ is equal to its relative double commutant in $C$, 
$A$ must be a von Neumann algebra.  A consequence
of this (and some other results from Sections 5 and 6) is  characterization of pairs $(A,B)$ admitting 
duality, deduced in Section 6.

In Section 5 we  introduce a suitable class of categories of operator bimodules and 
characterize cogenerators in such categories. The most important cases are of course
the categories of all and (if $A$ and $B$ are von Neumann algebras) of all normal operator 
$A,B$-bimodules.

It is natural to try to find the `smallest' possible among cogenerators. Let us call a strict cogenerator $Z_0$ 
{\em initial} if $Z_0$ is contained completely isometrically as an $A,B$-bimodule in every
strict cogenerator. If  all families of disjoint cyclic representations of $A$ and $B$ are 
countable, then the 
category of all  operator $A,B$-bimodules has an  initial  strict cogenerator (Section 7). 
This condition  is quite 
restrictive since separable C$^*$-algebras satisfying it turn out to be  a special subclass
of type I C$^*$-algebras. On the other hand, an analogous 
condition for a von Neumann algebra $M$, that the center of  $M$ is $\sigma$-finite
(each orthogonal family of nonzero central projections is countable), is not so restrictive. 
Define a cogenerator $Z_0$ in the category $\mnon$ of normal operator bimodules over  
von Neumann
algebras $M$ and $N$ to be {\em countably initial} if for each cogenerator $X$ there is
a completely isometric $M,N$-bimodule map from $Z_0$ into $X^{\bn}$ (= the $\ell^{\infty}$-direct
sum of countably many copies of $X$).  
If  the centers of $M$ and $N$ are $\sigma$-finite,  we shall specify
in Section 7 a countably initial cogenerator in $\mnon$. 
If the center of $M$ is not $\sigma$-finite, we shall show that the category ${{_M{\rm NOM}_{\bc}}}$
has no countably initial cogenerators.

The results here show that bimodules of the form $\bkh$
(with $\h$, $\k$ cyclic) are in some sense minimal among  (injective) cogenerators, and
the duality for such range bimodules is developed in \cite{Po}, \cite{M}. Apart from cogenerators,
bicommutation is the only aspect of duality studied in the present paper. 

Besides the  basic  operator space theory (which can be found in
the initial chapters of any of the books
\cite{BLM}, \cite{ER}, \cite{P}, \cite{Pi}),  we shall only need the notions
of an injective envelope and of a multiplier of an operator space, which we now
recall (more details are in  \cite[Chapter 4]{BLM} and \cite[Chapter 15]{P}).

\section{Preliminaries and notation}

Throughout the paper $A$ and $B$ are C$^*$-algebras.  $\aob$ denotes the category 
of all {\em operator 
$A,B$-bimodules} and  $\acbbb{X}{Y}$ the space of all completely bounded $A,B$-bimodule 
maps from $X$ to $Y$. The category of {\em Hilbert
$A$-modules} (that is, Hilbert spaces with the column operator space structure on which 
$A$ acts as a C$^*$-algebra) is 
denoted by $\ha$, while $\bm{A}{\k}{\h}$ means the space of all (completely) bounded 
$A$-module maps from $\k$ to $\h$. $\bkh$ and $\kkh$ are the spaces of all bounded and all compact
operators, respectively, from $\k$ to $\h$. For an operator $T\in\bkh$ and a cardinal
$\bi$ we denote by $T^{(\bi)}$ the direct sum of $\bi$ copies of $T$ and, for a subset
$\mathcal{S}\subseteq\bkh$, $\mathcal{S}^{(\bi)}$ is the set $\{T^{(\bi)}:\ T\in
\mathcal{S}\}$.
By $\rowi{X}$,  $\coli{X}$ and $\mati{X}$ we mean the set
of all rows, columns and $\bi\times\bi$ matrices, respectively, indexed by a set $\bi$, 
with the entries in an operator space $X$,
that represent bounded operators.

The {\em injective envelope} (\cite{H1},  \cite{Ru}, \cite{BLM}, \cite{ER},  
\cite{P}) of an operator space $X$ 
is an injective operator space $Y=I(X)$
containing $X$ such that $Y$ is the only injective subspace of $Y$ containing $X$.
An injective operator space $Y$ containing $X$ is the injective envelope of $X$
if and only if $Y$ has one (hence both) of the following two  
equivalent properties (see e.g. \cite[4.2.4]{BLM} or \cite[15.8]{P})
for a proof).

(i) {(\em Rigidity)} If $\phi:Y\to Y$ is a complete contraction such that $\phi|X$ is the
identity, then $\phi$ is the identity.

(ii) {(\em Essentiality)} If $\phi:Y\to Z$ is a complete contraction, where $Z$ is any
operator space, such that $\phi|X$ is a complete isometry, then $\phi$ is a complete isometry.

{\em A left multiplier of an operator space $X$} is a map $\phi:X\to X$ such that there
exist a Hilbert space $\h$, a complete isometry $\iota:X\to\bh$ and an
operator $T\in\bh$ such that $\iota\phi(x)=T\iota(x)$ for all $x\in X$. (That is, 
identifying $X$ with $\iota(X)$, $\phi$ is the restriction
to $X$ of a left multiplication by $T$). If, in addition, $T^*\iota(X)\subseteq\iota(X)$,
then $\phi$ is called an {\em adjointable left multiplier} on $X$. The multiplier norm 
$\norm{\phi}_m$ of $\phi$ is
the infimum of $\norm{T}$ over all possible such representations of $\phi$. 

In the theorem quoted below we recall two beautiful canonical descriptions of 
multipliers obtained by Blecher, Effros and Zarikian \cite{BEZ} and Blecher and Paulsen \cite{BP}
(see \cite[4.5.2]{BLM} or \cite[16.4]{P} for simpler proofs). 
Let 
$$S(X)=\left[\begin{array}{ll}
\bc&X\\
X^*&\bc
\end{array}\right]$$
be the operator system of $X$. Then the injective envelope $I(S(X))$ of $I(X)$ is
a C$^*$-algebra which can be decomposed as 
$$I(S(X))=\left[\begin{array}{cc}
I_{11}(X)&I(X)\\
I(X)^*&I_{22}(X)
\end{array}\right],$$
where $I(X)$ is the injective envelope of $X$ and $I_{11}(X)$ and $I_{22}(X)$ are
$C^*$-algebras over which $I(X)$ is a Hilbert C$^*$-bimodule.

\begin{theorem}\label{BEZ} A linear map $\phi$ on $X$ is a left multiplier
with $\norm{\phi}_m\leq1$ if and only if one of the following two equivalent conditions
holds:

(i) (\cite{BEZ}) The map
$$\tau:{\rm C}_2(X)\to {\rm C_2}(X),\ \ \ \tau\left(\left[\begin{array}{l}
x\\
y
\end{array}\right]\right)=\left[\begin{array}{c}
\phi(x)\\
y
\end{array}\right]$$
is completely contractive.

(ii) (\cite{BP}) There exists (a unique) $a\in I_{11}(X)$ with $\norm{a}\leq1$ such that $\phi(x)=ax$
for all $x\in X$.
\end{theorem}

The set $M_{\ell}(X)$ of all left multipliers of $X$ is an operator algebra and the subalgebra
$\al(X)$ of all adjointable multipliers is a C$^*$-algebra. 

Given Hilbert modules $\h\in\ha$, $\k\in\hb$ and $\l\in\ha\cap\hb$,
each completely bounded $A,B$-bimodule map $\phi:\kkh\to\B(\l)$ is of the form 

\begin{equation}\label{21} \phi(x)=U^*x^{(\bi)}V=\sumi u_i^*xv_i\ \ (x\in\kkh),
\end{equation}
where $\bi$ is an index set, $U\in\bm{A}{\l}{\h^{\bi}}$, $V\in\bm{B}{\l}{\k^{\bi}}$
are such that $\norm{U}\norm{V}=\cbnorm{\phi}$ and $u_i\in\bm{A}{\l}{\h}$ and
$v_i\in\bm{B}{\l}{\k}$ are the components of $U$ and $V$, respectively, when $U$
and $V$ are regarded as columns of operators using the identification 
$\bm{A}{\l}{\h^{\bi}}=\coli{\bm{A}{\l}{\h}}$. 

The proof of (\ref{21}) can be reduced to the special case $\k=\h$ by putting
$\g=\h\oplus\k$ (an $A\oplus B$-module), regarding $\kkh$ as the $(1,2)$ corner of $K:={\rm K}(\g)$ in the usual way
and extending $\phi$ by $0$ to the other three corners of $K$ to get a
bimodule map $\tilde{\phi}:K\to\B(\l)$.
Now we can refer to more general results (\cite[1.2]{M1}, \cite{S}) or reason as follows. 
By the representation theorem $\tilde{\phi}$
is of the form
$$\tilde{\phi}(x)=U^*\sigma(x)V \ \ \ (x\in K)$$
for a representation $\sigma$ of $K$ on a Hilbert space $\h_{\sigma}$ and operators
$U,V\in\bb{\l}{\h_{\sigma}}$ with $\norm{U}\norm{V}\leq\cbnorm{\tilde{\phi}}=\cbnorm{\phi}$,
whereby we may assume in addition that $$[\sigma(K)V\l]=\h_{\sigma}=[\sigma(K)U\l].$$ 
Since every representation
of $K$ is a multiple of the identity  \cite{Ar}, we
may assume that $\h_{\sigma}=\g^{\bi}$  and $\sigma(x)=x^{(\bi)}$
($x\in K$) for an index set $\bi$. Further, since $\tilde{\phi}$ is an $A$-module map, the equality
$[K^{(\bi)}V\l]=\g^{\bi}$ implies (by a standard computation)
that $U$ is an $A$-module map. 
Similarly $V\in\bm{B}{\l}{\g^{\bi}}$. Finally, with $P:\g\to\h$ and $Q:\g\to\k$
the orthogonal projections, we have  $x=PxQ$ if $x\in\kkh$. Hence, replacing $U$ and $V$ by $P^{(\bi)}U$ and $Q^{(\bi)}V$, respectively, 
yields $U\in\bm{A}{\l}{\h^{\bi}}$ and $V\in\bm{B}{\l}{\k^{\bi}}$
satisfying (\ref{21}).

\section{Embeddings of $\bkh$ into operator bimodules}

\begin{theorem}\label{th21} Let $\h\in\ha$ and $\k\in\hb$ be cyclic with unit cyclic
vectors $\xi\in\h$ and $\eta\in\k$ and let $\phi$ be a completely contractive
$A,B$-bimodule map from $\bkh$ (or $\kkh$) into an operator $A,B$-bimodule $X$. Then $\phi$ is completely isometric if and only if 
$\norm{\phi(\xi\otimes\eta^*)}=1$.
\end{theorem}

The Theorem is an immediate consequence of the following Lemma (in fact, of the identity
(\ref{29}), proved in the third paragraph of the proof of the Lemma). 
In later sections only the Theorem will be used, but the Lemma is a more complete result
with a consequence (Corollary \ref{co2}) that can not be deduced from the Theorem.

\begin{lemma}\label{le21} In the situation of Theorem \ref{th21}
there exists a Hilbert module $\l\in\ha\cap\hb$ such that:

(i) $X\subseteq\B(\l)$
completely isometrically as an operator bimodule. 

(ii) If $\phi|\kkh$, regarded as a map
into $\B(\l)$, is represented in the form (\ref{21}), where $U\in\bm{A}{\l}{\h^{\bi}}$
and $V\in\bm{B}{\l}{\k^{\bi}}$ are contractions, then $\norm{\phi(\xi\otimes\eta^*)}=1$
if and only if there exists an $\alpha\in\coli{\bc}$, with $\norm{\alpha}=1$, such that
\begin{equation}\label{22} UU^*\alpha=\alpha\ \ \ \mbox{and}\ \ \ VV^*\alpha=\alpha.
\end{equation}
In this case 
\begin{equation}\label{23}u:=U^*\alpha\in\bm{A}{\h}{\l}\ \ \ \mbox{and}\ \ \ v:=V^*\alpha
\in\bm{B}{\k}{\l}
\end{equation}
are isometries and
\begin{equation}\label{24} \phi(x)=uxv^*+\ort{e}\phi(x)\ort{f}\ \ (x\in\bkh),
\end{equation}
where $e=uu^*$ and $f=vv^*$. 
\end{lemma}

\begin{proof} If $\alpha$, as stated in the Lemma, does exist, then using (\ref{22}) and
the relations
$\norm{\alpha\xi}=1=\norm{\alpha\eta}$, $\norm{U}\leq1$ and $\norm{V}\leq1$, we have
for each $x\in\kkh$
$$\norm{\phi(x)}=\norm{U^*x^{(\bi)}V}\geq\inner{UU^*x^{(\bi)}VV^*\alpha\eta}{\alpha\xi}
=\inner{x^{(\bi)}\alpha\eta}{\alpha\xi}=\inner{x\eta}{\xi}.$$
This implies in particular that $\norm{\phi(\xi\otimes\eta^*)}=1$.

For the converse, first recall that by the CES representation theorem for operator bimodules
\cite[3.3.1]{BLM} there exists a C$^*$-algebra $D$ containing $X$ completely isometrically
as an operator $A,B$-bimodule (where the $A,B$-bimodule structure on $D$ is induced by
a pair of $*$-homomorphisms $A\to D$ and $B\to D$). Let $\l$ be the Hilbert space of the 
universal representation of $D$, regard
$D$ as contained in $\B(\l)$ and $\phi$ as a map into $\B(\l)$. Set $T=\phi(\xi\otimes\eta^*)$
and suppose that $\norm{T}=1$ Since each norm $1$ linear functional on $D$ is induced by
a pair of unit vectors in $\l$ (see \cite[10.1.3]{KR}), there exists unit vectors $\zeta,\tau
\in\l$ such that $\inner{T\tau}{\zeta}=\norm{T}=1$, hence from (\ref{21})
\begin{equation}\label{25} 1=\inner{\phi(\xi\otimes\eta^*)\tau}{\zeta}=
\inner{(\xi\otimes\eta^*)^{(\bi)}V\tau}{U\zeta}\leq\norm{V\tau}\norm{U\zeta}\leq1.
\end{equation}
Since equality occurs in the Cauchy - Schwarz inequality only for colinear vectors, 
(\ref{25}) implies that
\begin{equation}\label{26} (\xi\otimes\eta^*)^{(\bi)}V\tau=U\zeta\ \ \ \mbox{and similarly}
\ \ \ (\eta\otimes\xi^*)^{(\bi)}U\zeta=V\tau.
\end{equation}
Let $u_i$ and $v_i$ be the components of $U\in\coli{\bm{A}{\l}{\h}}$ and $V\in
\coli{\bm{B}{\l}{\k}}$, respectively, and set $\alpha_i=\inner{u_i\zeta}{\xi}$. Then,
written in components, the second identity in (\ref{26}) is
\begin{equation}\label{27} v_i\tau=\alpha_i\eta.
\end{equation}
Inserting (\ref{27}) in the first identity of (\ref{26}) shows that
\begin{equation}\label{28} u_i\zeta=\alpha_i\xi.
\end{equation}
From (\ref{25}) we also conclude that $\norm{U\zeta}=1=\norm{V\tau}$, which together with
(\ref{27}) (or (\ref{28})) implies that $\sumi|\alpha_i|^2=1$. Thus, the vector
$$\alpha:=(\alpha_i)\in\coli{\bc}$$
has norm $1$, so $\alpha\xi\in\h^{\bi}$ and $\alpha\eta\in\k^{\bi}$ are also unit vectors and
$U\zeta=\alpha\xi$, $V\tau=\alpha\eta$ by (\ref{28}) and (\ref{27}). Since $U$ and $V$ are
contractions, this implies that $U^*U\zeta=\zeta$ and $V^*V\tau=\tau$, which can be written
as $U^*(\alpha\xi)=\zeta$ and $V^*(\alpha\eta)=\tau$, hence
$$UU^*\alpha\xi=\alpha\xi\ \ \ \mbox{and}\ \ \ VV^*\alpha\eta=\alpha\eta.$$
Since $UU^*\alpha$ is an $A$-module map and $\xi$ is cyclic for $A$, the first of the above two identities
implies that $UU^*\alpha=\alpha$. Similarly $VV^*\alpha=\alpha$, which proves (\ref{22}).

With $u:=U^*\alpha$ and $v:=V^*\alpha$, we have from (\ref{22}) that $u^*u=\alpha^*UU^*\alpha
=\alpha^*\alpha=1$ and similarly $v^*v=1$. So $u$ and $v$ are isometries; let $e:=uu^*$ and
$f:=vv^*$ be their range projections. Using (\ref{21}), (\ref{23}) and (\ref{22}) we 
compute for $x\in\kkh$:
$$e\phi(x)f=uu^*U^*x^{(\bi)}Vvv^*=u\alpha^*UU^*x^{(\bi)}VV^*\alpha v^*=u\alpha^*x^{(\bi)}
\alpha v^*=\alpha^*\alpha uxv^*.$$
Thus
\begin{equation}\label{29} e\phi(x)f=uxv^*.
\end{equation}

To simplify the notation, we shall now regard $\h$ and $\k$ as subspaces of $\l$ using the
isometries $u$ and $v$. Then, relative to the decompositions $\l=\h\oplus\ort{\h}$ and
$\l=\k\oplus\ort{\k}$, the map $\phi$ can be represented as
$$\phi(x)=\left[\begin{array}{ll}
\phi_{11}(x)&\phi_{12}(x)\\
\phi_{21}(x)&\phi_{22}(x)
\end{array}\right]\ \ \ (x\in\bkh),$$
where $\phi_{ij}$ are complete contractions  and $\phi_{11}(x)=x$ for all $x\in\kkh$ by
(\ref{29}). Since $\bkh$ is the injective envelope of $\kkh$ \cite[4.6.12]{BLM}, $\phi_{11}$ must be the
identity by rigidity. 

To show that $\theta:=\phi_{21}=0$, note that 
$$\norm{x^*x+\theta(x)^*\theta(x)}\leq\norm{x^*x}\ \ \ (x\in\bkh)$$
since the first column of $\phi$ represents a contraction and $\phi_{11}={\rm id}$.
Choosing for $x$ any partial isometry and denoting by $p_x$ its initial projection $x^*x$,
it follows that $\norm{p_x+p_x\theta(x)^*\theta(x)p_x}\leq1$, hence $\theta(x)p_x=0$
and therefore
\begin{equation}\label{20}\theta(x)x^*=0.\end{equation}
If $\dim\k\leq\dim\h$, this implies in particular that $\theta(x)=0$ for each isometry
$x\in\bkh$, hence $\theta=0$ since the unit ball of $\bkh$ is the closure of the convex hull 
of isometries (by a polar decomposition argument and the fact that the open unit ball of a unital
C$^*$-algebra 
is the  convex hull of unitaries \cite[10.5.91]{KR}). If $\dim\k>\dim\h$, let $\k_1$ and 
$\k_2$ be any orthogonal subspaces of
$\k$ such that there exist partial isometries $x_i\in\bkh$ ($i=1,2$) with 
orthogonal ranges and the initial spaces equal to $\k_i$. Then $x_1+x_2$ is a partial
isometry (with the initial space $\k_1\oplus\k_2$), hence $\theta(x_1+x_2)(x_1+x_2)^*=0$
by (\ref{20}).
Using (\ref{20}) again, this implies $\theta(x_1)x_2^*+\theta(x_2)x_1^*=0$ and then, replacing
$x_2$ by $ix_2$, it follows that $\theta(x_2)x_1^*=0$, hence $\theta(x_2)|\k_1=0$. Since
$\ort{\k}_2$ is spanned by subspaces of the type $\k_1$, $\theta(x_2)|\ort{\k}_2=0$. But
the equality $\theta(x_2)x_2^*=0$ tells us that $\theta(x_2)|\k_2=0$, hence $\theta(x_2)=0$.
We conclude that $\theta(x)=0$ for each non-surjective partial isometry in $\bkh$. 
Since each coisometry is a sum of two such
partial isometries, this implies $\theta(x)=0$ for each coisometry $x\in\bkh$. Since $\bkh$ 
is the closure of the span of coisometries, $\phi_{21}=\theta$
must be $0$.
Similarly $\phi_{12}=0$, which proves that the above matrix of $\phi(x)$ is diagonal for
each $x\in\bkh$. 
Returning to the original notation, this proves the decomposition (\ref{24}) of $\phi$.
\end{proof}

\begin{remark}\label{re22} In the situation of Lemma \ref{le21}(ii) 
there exist two families $\{u_{\alpha}:\ \alpha\in\ba\}$ and $\{v_{\alpha}:\ 
\alpha\in\ba\}$ of isometries $u_{\alpha}\in\bm{A}{\h}{\l}$ and $v_{\alpha}\in\bm{B}{\k}{\l}$
such that:

(1) $u_{\alpha}^*u_{\beta}=0$ and $v_{\alpha}^*v_{\beta}=0$ if $\beta\ne\alpha$;

(2) $\phi(x)=\suma u_{\alpha}xv_{\alpha}^*+\ort{p}\phi(x)\ort{q}$, with $p:=\suma u_{\alpha}
u_{\alpha}^*$ and $q:=\suma v_{\alpha}v_{\alpha}^*$;

(3) $\norm{\ort{p}\phi(\xi_1\otimes\eta_1^*)\ort{q}\tau}<1$ for all unit cyclic vectors
$\xi_1\in\h$ and $\eta_1\in\k$ and all unit vectors $\tau\in\l$.

To prove this, let $\ba\subseteq\coli{\bc}$ be a maximal orthogonal family of vectors $\alpha$
satisfying the two equalities in (\ref{22}). For each $\alpha\in\ba$,  $u_{\alpha}:=U^*\alpha$
and $v_{\alpha}:=V^*\alpha$ are isometries by what we have already proved in Lemma \ref{le21}; 
let $e_{\alpha}=u_{\alpha}u_{\alpha}^*$ and $f_{\alpha}=
v_{\alpha}v_{\alpha}^*$ be their range projections. If $\beta\ne\alpha$ are in $\ba$ then,
using (\ref{22}) and the orthogonality of the set $\ba$,
$$u_{\beta}^*u_{\alpha}=\beta^*UU^*\alpha=\beta^*\alpha=0.$$
This shows that the family of projections $\{e_{\alpha}:\ \alpha\in\ba\}$ is orthogonal
and the same holds also for $\{f_{\alpha}:\ \alpha\in\ba\}$. Let
$$p=\suma e_{\alpha}\ \ \ \mbox{and}\ \ \ q=\suma f_{\alpha}.$$
Multiplying the identity $\phi(x)=u_{\alpha}xv_{\alpha}^*+
\ort{e_{\alpha}}\phi(x)\ort{f}_{\alpha}$ (which is just (\ref{24}) stated for each index
$\alpha$) from the left by $e_{\alpha}$ and from the right by $f_{\beta}$, where $\beta\ne\alpha$,
we get $e_{\alpha}\phi(x)f_{\beta}=0$. Similarly $e_{\alpha}\phi(x)\ort{q}=0$,
$\ort{p}\phi(x)f_{\alpha}=0$ and it follows that 
$$\phi(x)=\suma e_{\alpha}\phi(x)f_{\alpha}+\ort{p}\phi(x)\ort{q}=\suma u_{\alpha}x
v_{\alpha}^*+\ort{p}\phi(x)\ort{q}\ \ \ (x\in\bkh).$$

Suppose that there existed unit cyclic vectors $\xi_1\in\h$ and $\eta_1\in\k$ and a unit
vector $\tau\in\l$ such that
$\norm{\ort{p}\phi(\xi_1\otimes\eta_1^*)\ort{q}\tau}=1$. Then $\tau\in\ort{q}\l$ and,
denoting by $\zeta\in\ort{p}\l$ the unit vector such that $\inner{\ort{p}\phi(\xi_1\otimes\eta_1^*)
\ort{q}\tau}{\zeta}=1$,  we may apply the arguments
from the proof of Lemma \ref{le21} (from (\ref{25}) on) to the map 
$x\mapsto \ort{p}\phi(x)\ort{q}=(U\ort{p})^*x^{(\bi)}(V\ort{q})$. It follows that
there exists $\beta\in\coli{\bc}$ of norm $1$ such that
\begin{equation}\label{210} U\ort{p}U^*\beta=\beta\ \ \ \mbox{and}\ \ \ V\ort{q}V^*\beta=
\beta.
\end{equation}
The first of these two equalities implies that $\beta^*(1-U\ort{p}U^*)\beta=0$.
Since $0\leq1-UU^*\leq1-U\ort{p}U^*$, it follows that $\beta^*(1-UU^*)\beta=0$, hence
$(1-UU^*)\beta=0$ (by positivity) and $UU^*\beta=\beta$. Similarly $VV^*\beta=\beta$,
hence $\beta$ satisfies (\ref{22}). With $u_{\beta}$ and $v_{\beta}$ defined by
$$u_{\beta}=\ort{p}U^*\beta\ \ \ \mbox{and}\ \ \ v_{\beta}=\ort{q}V^*\beta,$$
it follows from (\ref{210}) that $u_{\beta}$ and $v_{\beta}$ are isometries 
(since $\norm{\beta}=1$). Moreover, since $u_{\alpha}=e_{\alpha}u_{\alpha}$ for each
$\alpha\in\ba$, we have
$$u_{\alpha}^*u_{\beta}=u_{\alpha}^*e_{\alpha}\ort{p}U^*\beta=0.$$
But, on the other hand, $u_{\alpha}^*u_{\beta}=\alpha^*U\ort{p}U^*\beta=\alpha^*\beta$
by (\ref{210}). It follows that $\alpha^*\beta=0$. This means that $\beta\perp\alpha$
for all $\alpha\in\ba$, which contradicts the maximality of the set $\ba$.
\end{remark}

Since each compact operator on a Hilbert space achieves its norm, the arguments
from Remark \ref{re22} and the proof of Lemma \ref{le21} imply the following corollary.

\begin{corollary}\label{co2} If $\h\in\ha$ and $\k\in\hb$ are cyclic with unit cyclic vectors
$\xi\in\h$ and $\eta\in\k$, and $\phi\in\acbbb{\kkh}{{\rm K}(\l)}$ is completely contractive 
(for an $\l\in\ha\cap\hb$) such 
that $\norm{\phi(\xi\otimes\eta^*)}
=1$, then for some (necessarily finite) $n$ there are orthogonal 
decompositions  $\l\cong\h^{n}\oplus\h_1$ (as $A$-modules)
and $\l\cong\k^{n}\oplus\k_1$ (as $B$-modules) relative to which $\phi$ is of the form
$\phi(x)=x^{(n)}\oplus\psi(x)$, where $\psi$ 
satisfies $\norm{\psi(\xi_1\otimes\eta_1^*)}<1$ for all unit cyclic vectors $\xi_1\in\h$
and $\eta_1\in\k$. 
\end{corollary}

\section{Induced maps on multiplier C$^*$-algebras}

In this section we show how complete contractions induce maps between the corresponding
multiplier C$^*$-algebras, but we shall only consider a special
situation needed later in the paper. 

\begin{theorem}\label{th31} Suppose that $\h\in\ha$, $\k$ is a Hilbert space, let $\eta\in\k$ 
be a unit vector,
$(\varepsilon_i)_{i\in\bi}$ an orthonormal basis of $\h$ and $x_i=\varepsilon_i\otimes\eta^*$.
Further, let  $X$ be an injective operator space and operator left $A$-module,  and let
$\phi:\kkh\to X$ be a completely contractive $A$-module map. Then the formula
\begin{equation}\label{32} \pl(T)=\sumi\phi(Tx_i)\phi(x_i)^*\ \ \ 
(T\in\kh)
\end{equation}
defines a completely contractive, completely positive $A$-bimodule map $\pl:\kh\to\al(X)$.
If $\phi$ is injective,  the kernel of $\pl$ does not
contain any positive (nonzero) elements. Moreover, if $\phi$ is completely 
isometric and $\h$ is cyclic, then $\pl$
is completely isometric. Finally, there exists an extension of $\pl$ to a map 
$\pl:\bh\to\al(X)$ with the same
properties.
\end{theorem}

\begin{proof} Since $\sumi x_ix_i^*=1_{\h}$ in the strong operator topology, for 
each $x\in\kkh$ we have that
\begin{equation}\label{311}x=\sumi x_ix_i^*x,
\end{equation}
where the sum is norm convergent.
For each finite subset $\bff$ of $\bi$ denote by $[x_i]_{\bff}$ the row matrix
in ${\rm R}_{\bff}(\kkh)$ with the entries $x_i$. Since $\phi$ is a complete contraction
and $\norm{[x_i]_{i\in\bi}}\leq1$, for each $T\in\kh$ we have
$$\begin{array}{lllll}\norm{\sum_{i\in\bff}\phi(Tx_i)\phi(x_i)^*}&=&
\norm{[\phi(Tx_i)]_{\bff}[\phi(x_i)]_{\bff}^*}
&\leq&\norm{[Tx_i]_{\bff}}\norm{[x_i]_{\bff}}\\
&\leq&\norm{\sum_{i\in\bff}Tx_i(Tx_i)^*}^{1/2}&
\leq&
\norm{T}.\end{array}$$
Since the sum $\sumi Tx_ix_i^*T^*$ is norm convergent, this estimate implies that also
the sum (\ref{32})
is  norm convergent and defines a contraction $\pl:\kh\to\kl(X):=[XX^*]$. 
A similar argument (available also in a more general context of Remark \ref{ret}(ii)
below) shows that $\pl$ is completely contractive. But here an even simpler
argument is possible: if $[t_{ij}]$ is the matrix of $T$ relative to the orthonormal basis 
$(\varepsilon_i)_{i\in\bi}$, then by a short computation $\pl(T)$ is just a product
of three operator matrices
\begin{equation}\label{312} \pl(T)=
[\phi(x_i)][t_{ij}][\phi(x_j)]^*,
\end{equation}
from which one can see that $\pl$ is completely contractive and completely positive. 
Since $X$ is injective, it is known  that $[XX^*]$ is a C$^*$-subalgebra of
$\al(X)$ and that $\al(X)$ is injective. (See Theorem \ref{BEZ}(ii) and the paragraph
above Theorem \ref{BEZ} or \cite[4.4.3]{BLM}.  $\al(X)$ is known
to be just the multiplier C$^*$-algebra of $[XX^*]$.)

Since $\phi$ is a left $A$-module map, $\pl$ is clearly a left $A$-module map by (\ref{32}). To
show that $\phi$ is also a right $A$-module map, we use the specific form of the operators
$x_i$. For each $T\in\kh$ of the form $T=\xi\otimes\zeta^*$ ($\xi,\zeta\in\h$) we have
$$\sumi\phi(Tx_i)\phi(x_i)^*=\sumi\phi(\inner{\varepsilon_i}{\zeta}\xi\otimes\eta^*)
\phi(\varepsilon_i\otimes\eta^*)^*=
\phi(\xi\otimes\eta^*)\phi(\sumi\inner{\zeta}{\varepsilon_i}\varepsilon_i\otimes\eta^*)^*,$$
hence, by (\ref{32}) and since $(\varepsilon_i)_{i\in\bi}$ is an orthonormal basis of $\h$,
\begin{equation}\label{33} \pl(\xi\otimes\zeta^*)=\phi(\xi\otimes\eta^*)\phi(\zeta\otimes
\eta^*)^*.
\end{equation}
Therefore, for $T$ of the form $\xi\otimes\zeta^*$ and $a\in A$ we have
$$\pl(Ta)=\pl(\xi\otimes(a^*\zeta)^*)=\phi(\xi\otimes\eta^*)\phi(a^*\zeta\otimes\eta^*)^*=
\phi(\xi\otimes\eta^*)\phi(\zeta\otimes\eta^*)^*a=\pl(T)a$$
since $\phi$ is a left $A$-module map. Since operators of the form $\xi\otimes\zeta^*$
densely span $\kh$, this shows that $\pl$ is a right $A$-module map.  If $\phi$ is 
completely isometric, then
it follows from (\ref{33}) that $\norm{\pl(\xi\otimes\xi^*)}=1$ for all unit vectors
$\xi\in\h$, hence, if $\h$ is cyclic, $\pl$ is completely isometric by 
Theorem \ref{th21}. If $\phi$ is injective, then (\ref{33}) with $\zeta=\xi$ (together
with positivity) implies that the kernel of
$\pl$ does not contain any nonzero positive operator.

We may assume that $\al(X)\subseteq\B(\l)$ for some Hilbert space $\l$. Since $\al(X)$ is 
injective,  there exists a  completely contractive projection $E:\B(\l)\to\al(X)$, which is
automatically an $\al(X)$-bimodule map by a well known result of Tomiyama \cite[10.5.86]{KR}.
We can now extend $\pl$ to a (normal) completely positive $A$-bimodule complete contraction 
from $\bh$ to $\B(\l)$ and then compose this extension with $E$ to get the completely positive
complete contraction $\pl:\bh\to\al(X)$ with required properties. 
\end{proof}

\begin{remark}\label{ret} (i) Let $\h$ be a Hilbert space, $X$ a ternary ring of operators 
(see \cite{BLM}) and $\psi:\h\to X$ a completely bounded map. 
Denote by $\psi^*:\h^*\to X^*$ the map $\psi^*(\xi^*)=\psi(\xi)^*$ and by
$\tau_{\psi}:\kh\to\al(X)$ the composition of the maps in the diagram
$$\kh=\h\stackrel{h}{\otimes}\h^*\stackrel{\psi\otimes\psi^*}{\longrightarrow}X\stackrel{h}{\otimes}
X^*\stackrel{\mu}{\longrightarrow}[XX^*]\subseteq\al(X),$$
where $\mu$ is the multiplication: $\mu(x\otimes y^*)=xy^*$. For a Hilbert space $\k$
and a fixed unit vector $\eta\in\k$ let $\iota_{\eta}:\h\to\h\stackrel{h}{\otimes}\k^*$
be the embedding $\iota_{\eta}(\xi)=\xi\otimes\eta^*$. Given a  completely
bounded map $\phi:\kkh=\h\stackrel{h}{\otimes}\k^*\to X$, let  $\psi_{\eta}:\h\to X$ 
be $\psi_{\eta}=\phi\iota_{\eta}$. Then it can be verified that
$\tau_{\psi_{\eta}}$
is just the map $\pl$ defined by (\ref{32}) in Theorem \ref{th31}, which explains some of the
properties of $\pl$ and shows that $\pl$ is independent of the choice of orthonormal basis.
It depends, however, on the choice of $\eta$.

(ii) Using formula (\ref{32}) one can define a map in the situation when $\kkh$ is replaced
by a ternary ring of operators $W$  such that $W$ has a left frame. 
Here {\em a left frame} in  $W$ is a row $[x_i]\in\rowi{W}$ (for some index set
$\bi$) with $\norm{[x_i]}\leq1$ such that
$w=\sumi x_ix_i^* w$ for all $w\in W,$
where the sum is norm convergent.  
It follows from a result of Brown (\cite[2.2]{Br},  \cite[5.53]{RW}) that $W$ has a 
left frame if $[WW^*]$ has a countable approximate identity.
We do not know if $\pl$  can be constructed in a
more general situation, 
without assuming the existence of a left frame in $W$. 
\end{remark}

\section{On  cogenerators}

In the definition below we extract some properties which are common to certain natural 
categories of operator bimodules, such as the category $\aob$ of all  operator
$A,B$-bimodules and, if $M$ and $N$ are von Neumann algebras, the category $\mnon$
of all normal operator $M,N$-bimodules (as well as some other categories defined in \cite{M}).
The main feature of such categories is that they have enough bimodules of the form $\bkh$.
A reason for introducing them in the definition below is that their cogenerators are suitable for studying
a bicommutation problem (Theorem \ref{th43} below). But, if not interested in general 
categories of operator bimodules, the reader may skip a few paragraphs and continue reading at 
Definition \ref{dc}, with the concrete categories $\aob$ or $\mnon$ in mind.

\begin{remark}\label{re40} In (vii) of the following definition we will use the
elementary fact (which may be proved by using matrix units) that each Hilbert 
$\matn{A}$-module $\h$ is (up to a unitary equivalence) of the form  $\g^n$ for a
Hilbert $A$-module $\g$, and that an $\matn{A},\matn{B}$-bimodule map $\psi$ from
$\matn{X}$ into $\B(\l^n,\g^n)=\matn{\B(\l,\g)}$ is necessarily of the form $[x_{ij}]
\mapsto [\phi(x_{ij}]$ for a (unique) $A,B$-bimodule map $\phi:X\to\B(\l,\g)$.
We shall write this as $\psi=\phi_n$.
\end{remark}

\begin{definition}\label{d41} A  subcategory $\mathcal{O}$ of $\aob$ is called an {\em ample
category}
if it satisfies the following conditions:

(i) If $a\in A$ and $aX=0$ for all $X\in\om$ then $a=0$ and similarly for $b\in B$.

(ii) $\mati{X}\in\om$ for each $X\in\om$ and each index set $\bi$. 

(iii) For each index sets $\bi$, $\bj$
and all $a\in{\rm M}_{\bj,\bi}(\bc)$, $b\in{\rm M}_{\bi,\bj}(\bc)$  the map $\mu_{a,b}$
defined by $\mu_{a,b}(x)=axb$ is a morphism in $\om(\mati{X},
\matj{X})$.

(iv) Given a map $\phi=[\phi_{ij}]\in\cbb{X}{\mati{Y}}$ such that $\phi_{ij}\in\om(X,Y)$ for
all $i,j\in\bi$, we have that $\phi\in\om(X,\mati{Y})$.

(v) For any collection 
$(X_{i})_{i\in \bi}$ of modules in $\om$ the $\ell^{\infty}$-direct sum
$X:=\oplus_{i\in\bi} X_i$ is in $\om$ and the coordinate inclusions $X_i\to X$
and projections $X\to X_i$ are morphisms in $\om$. Moreover, for each $Y\in\om$ and bounded family
of morphisms $\phi_i\in\om(Y,X_i)$ the map $\phi$ defined by $\phi(y)=\oplus_{i\in\bi}
\phi_i(y)$ is in $\om(Y,X)$.

(vi) If $\bkh\in\om$ for some Hilbert modules $\h\in\ha$ and $\k\in\hb$ and if 
$\h_0$ and $\k_0$ are (isometrically) isomorphic to Hilbert submodules of $\h$ and $\k$, 
respectively,
then $\B(\k_0,\h_0)\in\om$ and the induced `inclusion' $\B(\k_0,\h_0)\to\B(\k,\h)$ 
and the compression $\B(\k,\h)\to\B(\k_0,\h_0)$ are morphisms in $\om$.

(vii) For each $X\in\om$, $n\in\bn$, $x\in\matn{X}$ and $\varepsilon>0$ there exist Hilbert modules 
$\h^n\in{{_{\matn{A}}{\rm HM}}}$ and $\k^n\in{{_{\matn{B}}{\rm HM}}}$
such that $\bkh\in \mathcal{O}$,
and there exists a completely contractive 
$\matn{A},\matn{B}$-bimodule map $\phi_n:\matn{X}\to\B(\k^n,\h^n)$ such that
$\phi\in\om(X,\bkh)$ and
$\norm{\phi_n(x)}>\norm{x}-\varepsilon$.
\end{definition}

The conditions (vii) and (v) together imply  that each bimodule $X$ in an ample
category can be completely isometrically embedded into an $\ell^{\infty}$-direct
sum of bimodules of the form $\bkh$ with a morphism within the category.

\begin{remark}\label{re4} (i) Using for $a$ and $b$ suitable matrices (with the
entries $1$ and $0$ only), \ref{d41}(iii) implies in particular
that if $\bj\subseteq\bi$,
the canonical inclusion $\matj{X}\to\mati{X}$  and the compression $\mati{X}\to\matj{X}$
are morphisms in $\om$. 

(ii) Since a composition of morphisms in a category is
a morphism, it follows from (iii) and (iv) of Definition \ref{d41} that $\om(X,Y)$ is a
vector subspace of $\acbbb{X}{Y}$. To see,
for example, that the sum $\phi+\psi$ of two morphisms in $\om(X,Y)$ is in $\om(X,Y)$,
note that the diagonal matrix $\delta:=\phi\oplus\psi$ 
is in $\om(X,{\rm M}_2(Y))$ by \ref{d41}(iv). Then, with $\mu_{a,b}
\in\om({\rm M}_2(Y),Y)$ defined as in \ref{d41}(iii), where $a=[1,1]$
and $b=[1,1]^T$, we have that $\phi+\psi=\mu_{a,b}\circ\delta$ is in $\om(X,Y)$.

(iii) For each $\phi\in\om(X,Y)$ the amplification $\phi_{\bi}:[x_{ij}]\mapsto[\phi(x_{ij})]$
is a morphism in $\om(\mati{X},\mati{Y})$.

To show this, let $E_i\in\coli{\bc}$ have $1$ on the $i$-th position and $0$ elsewhere
and let $E_i^T$ be its transposition. Then the map $\phi_{ij}(y):=\phi(E_i^TyE_j)$ is in
$\om(\mati{X},X)$ by  \ref{d41}(iii). Since the map $[\phi_{ij}]:\mati{X}\to\mati{X}$
is just the amplification
$\phi_{\bi}$ of $\phi$, it follows from  \ref{d41}(iv)
that $\phi_{\bi}\in\om(\mati{X})$.

(iv) If $\{\psi_i:\ i\in\bi\}$ is a bounded set of morphisms $\psi_i\in\om(Y_i,X_i)$,
then $\psi:=\oplus_{i\in\bi}\psi_i$ is a morphism in $\om(Y,X)$, where
$X=\oplus_{i\in\bi}X_i$ and $Y=\oplus_{i\in\bi}Y_i$. This follows by
applying \ref{d41}(v) to morphisms $\phi_i:=\psi_i\eta_i$, where $\eta_i\in\om(Y,Y_i)$
is the coordinate projection.
\end{remark}

{\em Injectivity} in  an ample category is defined as usual,
considering only completely contractive morphisms in the category.

\begin{remark}\label{re41} If $\om$ is ample and $X$ is injective in $\om$, then 
$X$ is injective as an operator space, hence injective in $\aob$. 

To see this, let $\iota$ be a completely isometric embedding of $X$ into  
an $\ell^{\infty}$-direct sum $Y\in\om$
of bimodules of the form $\bkh\in\om$ such that $\iota\in\om(X,Y)$  
(\ref{d41}(vii) and (v)).   
By  injectivity of $X$ in $\om$ there exists a completely contractive extension of the identity 
$1_X$ to a morphism $E\in\om(Y,X)$ (that is, $E\iota=1_X$). 
Hence $X$ must be injective as an operator space since $Y$ is.
\end{remark}

\begin{definition}\label{dc} If $\mathcal{O}$ is a  subcategory in $\aob$, a bimodule $Z\in\mathcal{O}$
is called a {\em cogenerator} if for each $Y\in\om$, $y\in Y$ and $\varepsilon>0$ there
exists a completely contractive morphism $\phi\in\om(Y,Z)$  such that 
$\norm{\phi(y)}>\norm{y}-\varepsilon$. If (in addition) $\phi$ can be chosen so that $\norm{\phi(y)}=
\norm{y}$, then $Z$ is called a {\em strict cogenerator}.
\end{definition}

For example, if $\h\in\ha$ and $\k\in\hb$ contain (up to a unitary equivalence) all cyclic Hilbert modules over 
$A$ and $B$, respectively, then $\bkh$ is an injective strict cogenerator in $\aob$. (That it
is a cogenerator can be deduced from the well known CES  theorem \cite[3.3.1]{BLM},
or from the operator bipolar theorem \cite[1.1]{M2}). That such a
$\bkh$ is in fact a strict cogenerator follows from \cite[4.1]{Po}, but we can also present an 
alternative argument.
Given $y$ in an operator bimodule $Y$, first choose a norm one functional $\rho$ on $Y$ such 
that
$\rho(y)=\norm{y}$. Then apply the well known CSPS factorization theorem to the map
$\tilde{\rho}:A\times Y\times B\to\bc$, $\tilde{\rho}(a,x,b):=\rho(axb)$, to show that
$\rho$ is of the form $\rho(axb)=\inner{\pi(a)\phi(x)\sigma(b)\eta}{\xi}$ for  some cyclic
representations $\pi:A\to\B(\h_{\pi})$ and $\sigma:B\to\B(\k_{\sigma})$, with unit cyclic
vectors $\xi$ and $\eta$, respectively, and a completely 
contractive $A,B$-bimodule map $\phi:Y\to\B(\k_{\sigma},\h_{\pi})$. Then 
$\inner{\phi(y)\eta}{\xi}=\rho(y)=\norm{y}$ implies that $\norm{\phi(y)}=\norm{y}$.
By cyclicity we can regard
$\h_{\pi}$ and $\k_{\sigma}$ as contained in $\h$ and $\k$, respectively. 

Similarly, if $A$ and $B$ are von Neumann algebras and
$\h$ and $\k$ are  Hilbert spaces of the universal normal representations of $A$ and
$B$, respectively, then $\bkh$ is an injective cogenerator in the category of all 
normal operator $A,B$-bimodules (and normal dual operator $A,B$-bimodules).

\begin{example} As an example of a cogenerator in $\aob$ that is not injective and not strict,
consider $A=B=c_0$ acting on $\ell^2$ in the usual way and let $Z={\rm K}(\ell^2)$.
To see that $Z$ is a cogenerator, first note that each cyclic representation of $c_0$ is
contained in the representation of $c_0$ on $\ell^2$ as diagonal matrices (namely, $\ell^{\infty}$
is the universal von Neumann envelope of $c_0$). This implies
that $\B(\ell^2)$ is a (strict) cogenerator. Now, let $p_n$ be the projection onto the first
$n$ coordinates in $\ell^2$. Since $p_n\in c_0\subseteq \com{c_0}$, the compression
$x\mapsto p_nxp_n$ ($x\in\B(\ell^2)$) is a $c_0$-bimodule map and the map
$x\mapsto\oplus_{n=1}^{\infty}p_nxp_n$ is a completely isometric $c_0$-bimodule map
from $\B(\ell^2)$ into ${\rm K}(\ell^2)^{\infty}$. This implies that ${\rm K}(\ell^2)$
is a cogenerator for operator $c_0$-bimodules. 

To show that ${\rm K}(\ell^2)$ is not a strict cogenerator, first note that
each $c_0$-bimodule complete
contraction $\phi:\B(\ell^2)\to{\rm K}(\ell^2)$ is weak* continuous (when considered as
a map into $\B(\ell^2)$), hence  of the form $\phi(x)=U^*x^{(\bi)}V$
for suitable contractions $U$ and $V$ (as in (\ref{21})). This follows from the decomposition
$\phi=\phi_n+\phi_s$ of $\phi$ into the normal and the singular part (which can be verified to
be $c_0$-bimodule maps) and noting that $\phi_s$
must be zero by the bimodule property. Namely, if $e_i$ are the minimal projections in $c_0$,
then (since $\phi_s$ annihilates compact operators; see \cite[Chapter 10,
10.5.15]{KR}) for each $x\in\B(\ell^2)$ we have
$0=\phi_s(e_ixe_j)=e_i\phi_s(x)e_j$, hence $\phi_s(x)=0$ since $\sum e_i=1_{\ell^2}$. 
Now choose an operator $y\in\B(\ell^2)$ which does not achieve its norm. Since $\phi(y)$ is compact, it achieves its norm
on a unit vector $\eta\in\ell^2$, so that $\norm{\phi(y)}=\norm{\phi(y)\eta}=
\norm{U^*y^{(\bi)}V\eta}\leq\norm{y}$. If $\norm{\phi(y)}=\norm{y}$, then it follows that
$\norm{U^*y^{(\bi)}V\eta}=\norm{y}$, hence $\norm{y^{(\bi)}V\eta}=\norm{y}$ or
\begin{equation}\label{4000}\sumi\norm{yv_i\eta}^2=\norm{y}^2,
\end{equation}
where $v_i$ are the components of $V$. But, since $\norm{V}\leq1$, 
$$\sumi\norm{yv_i\eta}^2
\leq\norm{y}^2\sumi\norm{v_i\eta}^2\leq\norm{y}^2$$
and therefore the equality (\ref{4000})
is possible only if $\norm{yv_i\eta}=\norm{y}\norm{v_i\eta}$, which means that $y$
achieves its norm at the unit vector $\norm{v_i\eta}^{-1}v_i\eta$ for some $i$. 
But this contradicts the choice of $y$.
\end{example}

\begin{proposition}\label{pr41} A bimodule $Z$ in an ample category $\om$ is a cogenerator
(is a strict cogenerator in $\om=\aob$) if and only if for each $\bkh\in\om$, with $\h\in\ha$ and $\k\in\hb$ 
cyclic, there exists
a complete isometry in $\om(\bkh,\zn)$ (a complete isometry in ${{\rm CB}_A(\bkh,Z)_B}$, 
respectively).
\end{proposition} 

\begin{proof} If $Z$ is a cogenerator in $\om$ and $\h\in\ha$, $\k\in\hb$ are cyclic,
with  unit cyclic vectors $\xi\in\h$ and $\eta\in\k$ and
such that $\bkh\in\om$, then for each $n=1,2,\ldots$ there exists a complete contraction
$\phi_n\in\om(\bkh,Z)$ such that $\norm{\phi_n(\xi\otimes\eta^*)}>1-
\frac{1}{n}$. Then by Theorem \ref{th21} the direct sum $\phi$ of the maps $\phi_n$
embeds $\bkh$ into $\zn$ completely isometrically and $\phi\in\om(\bkh,\zn)$ by
\ref{d41}(v). 

For the converse, let $Y\in\om$, $y\in Y$, with $\norm{y}=1$, and $\varepsilon>0$.
By  \ref{d41}(vii) there exist $\bkh\in\om$ and a
complete contraction $\phi\in\om(Y,\bkh)$  such that $\norm{\phi(y)}>1-\varepsilon$.
Let $\xi\in\h$ and $\eta\in\k$ be unit vectors such that $\inner{\phi(y)\eta}{\xi}>
1-\varepsilon$,  let $\h_0:=[A\xi]$ and $\k_0:=[B\eta]$ be the corresponding 
cyclic submodules and $p:\h\to\h_0$  the orthogonal projection. 
Then $\B(\k_0,\h_0)\in\om$, the
map $\phi_0(x):=p\phi(x)|\k_0$ from $Y$ into $\B(\k_0,\h_0)$ is a morphism in $\om$ by 
\ref{d41}(vi) and satisfies 
$\norm{\phi_0(y)}>1-\varepsilon$. By hypothesis there exists a complete isometry
$\psi\in\om(\B(\k_0,\h_0),\zn)$. Let $\theta=\psi\phi_0$, so that $\theta\in\om(Y,\zn)$
is a complete contraction.
Composing $\theta$ with a suitable coordinate projection $\zn\to Z$ (which is in $\om(\zn,Z)$
by  \ref{d41}(v)), we find a map
in $\om(Y,Z)$ satisfying the requirement in the definition of a cogenerator.

The proof for strict cogenerators in $\aob$ is similar (see the comment following 
Definition \ref{dc}).
\end{proof}

\begin{corollary} If $Z$ in an injective strict cogenerator in $\aob$, then $\al(Z)$ is an
injective strict cogenerator in ${_A{\rm OM}_A}$.
\end{corollary}

\begin{proof} Let $\h$ be any cyclic Hilbert module over $A$ and choose any such module $\k$
over $B$. By Proposition \ref{pr41} there exists a completely isometric $A,B$-bimodule
map $\phi:\bkh\to Z$, which by Theorem \ref{th31} induces a completely isometric $A$-bimodule
map $\pl:\bh\to\al(Z)$. By Proposition \ref{pr41} again this implies that $\al(Z)$ is
a strict cogenerator in ${_A{\rm OM}_A}$. It is known that $\al(Z)$ is injective if $Z$ is injective.
\end{proof}

\begin{theorem}\label{th4} If $Z$ is a cogenerator in an ample category $\om$, then for each
$X\in\om$ there exist an index set $\bk$ and a complete isometry in $\om(X,\matk{Z})$.
\end{theorem}

\begin{proof} From \ref{d41}(vii),(v) there exist two families of Hilbert modules $\h_m\in\ha$ and $\k_m\in\hb$
($m\in\bM$), such that $\B(\k_m,\h_m)\in\om$, $Y:=\oplus_{m\in\bM}\B(\k_m,\h_m)\in\om$,
and there exists a complete isometry $\iota\in\om(X,Y)$. Suppose that for each $m$ we have found
a complete isometry $\phi_m\in\om(\B(\k_m,\h_m),{\rm M}_{\bk_m}(Z))$ for some $\bk_m$.
Then, let $\bk$ be the disjoint union of the sets $\bk_m$, $W=\oplus_{m\in\bM}{\rm M}_{\bk_m}(Z)$,
$\eta_m:W
\to{\rm M}_{\bk_m}(Z)$ the coordinate projection (a morphism in $\om$ by  
\ref{d41}(v)), $\iota_m:{\rm M}_{\bk_m}(Z)\to {\rm M}_{\bk}(Z)$ the canonical inclusion and
$\tilde{\eta}_m:=\iota_m\eta_m$. Since $\iota_m$ is a morphism in $\om$ by Remark \ref{re4}(i),
the same holds for $\tilde{\eta}_m$ and then it follows from  \ref{d41}(iv) that
the diagonal matrix $\oplus_m\tilde{\eta}_m\in\cbb{W}{{\rm M}_{\bM}({\rm M}_{\bk}(Z))}$ 
represents
a morphism $\delta$ in $\om(W,{\rm M}_{\bM}({\rm M}_{\bk}(Z)))$. For each $k\in\bk$ denote
by $m_k$ the (unique) element of $\bM$ such that $k\in\bk_{m_k}$ and regard $\bk$ as a
subset in $\bM\times\bk$ by $k\mapsto(m_k,k)$. The compression 
$$\gamma:{\rm M}_{\bM}({\rm M}_{\bk}(Z))={\rm M}_{\bM\times\bk}(Z)\to{\rm M}_{\bk}(Z)$$
is a morphism in $\om$ by Remark \ref{re4}(i), hence so is $\gamma\delta$. But this map 
$\gamma\delta$ is just the canonical inclusion $\kappa$ of $W$ into (the block-diagonal
matrices in) ${\rm M}_{\bk}(Z)$.
Then the composition $\kappa\circ(\oplus_{m\in\bM}\phi_m)\iota$ is a morphism in $\om$
(Remark \ref{re4}(iv)) and
embeds $X$ completely isometrically into ${\rm M}_{\bk}(Z)$.
Thus, it suffices to find the appropriate maps $\phi_m$, which reduces the proof to showing that
each bimodule of the form $\bkh\in\om$ can be embedded completely isometrically into  
$\matk{Z}$ (for some $\bk$) by a morphism in $\om$.  

Denote by $\pi:A\to\bh$ and $\sigma:B\to\B(\k)$ the  representations that induce the
module structures on $\h$ and $\k$.
Let $\{\h_i:\ i\in\bi\}$ be a maximal
set of disjoint cyclic Hilbert submodules in $\h$. Here `disjoint' means that $\h_i$ and
$\h_j$ have no isomorphic non-zero submodules if $i\ne j$. Denoting by $\com{e}_i:\h
\to\h_i$ the orthogonal projections (thus $\com{e}_i\in\com{\pi(A)}$), disjointness 
means that the central carriers $p_i\in\overline{\pi(A)}$ of 
projections $\com{e}_i$ are mutually orthogonal \cite[10.3.3]{KR}. Let 
$\com{e}=\sumi\com{e}_i$. Since the central carrier of $\com{e}$ is $1_{\h}$ by maximality, the
map $\overline{\pi(A)}\to\overline{\pi(A)}\com{e}$, $a\mapsto a\com{e}$ is a $*$-isomorphism
of von Neumann algebras by \cite[5.5.5]{KR}, hence $\tilde{\h}:=\oplus_{i\in\bi}\h_i$ is a 
faithful $\overline{\pi(A)}$-module.
Therefore $\h$ is (isometrically isomorphic to) a submodule in $\tilde{\h}^{\bl}$ for
some $\bl$ (as a module over $\overline{\pi(A)}$, hence also over $A$). Denoting by $\tilde{\k}$ a submodule of $\k$ that is constructed in the same 
way as $\tilde{\h}$ (and enlarging $\bl$ if necessary), it follows that $\bkh$ is `contained'
in $\B(\tilde{\k}^{\bl},
\tilde{\h}^{\bl})={\rm M}_{\bl}(\B(\tilde{\k},\tilde{\h}))$, where the `inclusion' is a
morphism in $\om$ by  \ref{d41}(vi). Hence (simplifying the notation),
the proof reduces to  modules of the form $\bkh$, where $\h=\oplus_{i\in\bi}
\h_i$ and $\k=\oplus_{j\in\bj}\k_j$ are now direct sums of disjoint cyclic submodules. 

Let $\com{p}_i:\h\to\h_i$ and $\com{q}_j:\k\to\k_j$ be the projections. Let $\bs$ be
the set of all pairs $(\g,\l)$, where $\g\subseteq\h$ and $\l\subseteq\k$ are cyclic
Hilbert modules and for each $s=(\g,\l)\in\bs$ let $\psi_s$ be the compression $x\mapsto Px|\l$,
where $P:\h\to\g$ is the orthogonal projection. Note that $\psi_s\in\om(\bkh,\B(\l,\g))$ 
by  \ref{d41}(vi). By Proposition \ref{pr41} for
each $s=(\g,\l)\in\bs$ there is a complete isometry  $\phi_s\in\om(\B(\l,\g),\zn)$. Define 
$$\phi:\bkh\to(\zn)^{\bs},\ \ \ \phi(x)=\oplus_{s\in\bs}\phi_s(\psi_s(x))\ \ \ (x\in\bkh).$$
Then $\phi\in\om(\bkh,(\zn)^{\bs})$ by  \ref{d41}(v). To show that $\phi$ is 
isometric, note
(by finite rank approximation) that for a compact $x$ there exist
countable subsets $\bi_x\subseteq\bi$ and $\bj_x\subseteq\bj$ such that $\com{p}_ix=0$
if $i\notin\bi_x$ and $x\com{q}_j=0$ if $j\notin\bj_x$. By disjointness 
the modules $\h_x:=\oplus_{i\in\bi_x}\h_i$ and $\k_x:=\oplus_{j\in\bj_x}\k_j$ are
cyclic \cite[5.5.10]{KR}, so $s:=(\k_x,\h_x)\in\bs$. For this particular $s$ we have
that $\norm{\phi_s(\psi_s(x))}=\norm{x}$. This shows that $\phi|\kkh$ is isometric and a similar
argument shows that $\phi|\kkh$ is completely isometric. Since $\bkh$ is the injective
envelope of $\kkh$, $\phi$ must be a complete isometry.
\end{proof}

\section{On relative bicommutants}

We shall now turn to the characterization of pairs $(A,B)$ admitting duality, but first we 
need a general result concerning bicommutants and two short lemmas.

\begin{theorem}\label{th32} Let  $\pi:A\to C$ be a $*$-homomorphism between C$^*$-algebras
(making $C$ an $A$-bimodule). Assume that $C$ is injective  and  that $A\subseteq\bh$ for
some Hilbert space $\h$ such that there exists a bounded left $A$-module map 
$\phi:\kh\to C$  the kernel of which does not contain any
nontrivial right ideals. If $\pi(A)$ is equal to its relative double commutant
$\pi(A)^{cc}$ in $C$, then $A$ is a von Neumann algebra.
\end{theorem}

\begin{proof} We may assume that $C\subseteq\B(\l)$ for some Hilbert space $\l$ and denote
by $\iota$ the inclusion of $C$ into $\B(\l)$ and by
$E$ the conditional expectation from $\B(\l)$ to $C$. Let $\tilde{A}$ be the universal
von Neumann envelope of $A$,  $\Phi:A\to\tilde{A}$ the universal representation and let 
$\alpha:\tilde{A}\to\overline{A}$ be the weak* continuous extension of $\Phi^{-1}$
(see \cite[Section 10.1]{KR} if necessary). Let $\tilde{\pi}:\tilde{A}\to\B(\l)$ be
the weak* continuous extension of $\iota\pi\Phi^{-1}$ from $\Phi(A)$ to  $\tilde{A}$.
Since $\tilde{\pi}(T)\in\overline{\pi(A)}$ (= the weak* closure of $\pi(A)$ in $\B(\l)$)
and $E$ is a $C$-bimodule map, for each $c\in \pi(A)^c$ we have
$$E(\tilde{\pi}(T))c=E(\tilde{\pi}(T)c)=E(c\tilde{\pi}(T))=cE(\tilde{\pi}(T)),$$
hence $E(\tilde{\pi}(T))\in\pi(A)^{cc}$. By hypothesis there exists an $a_T\in A$
such that
\begin{equation}\label{34}E(\tilde{\pi}(T))=\pi(a_T).
\end{equation}

For each $x\in\kh$ the map 
$\phi_x:\tilde{A}\to\B(\l)$, $\phi_x(T):=\iota\phi(\alpha(T)x)$,
is weak* continuous, since $\phi_x$ is the composition of $\alpha$, 
the right multiplication by
$x$ and $\iota\phi$. (Note that $\iota\phi$, as any bounded linear map 
on $\kh$, is weak* continuous
since this holds for bounded linear functionals). Further, the map $\psi_x:\tilde{A}\to
\B(\l)$, $\psi_x(T)=\tilde{\pi}(T)\iota\phi(x)$ is also weak* continuous since this holds for
$\tilde{\pi}$. If $T=\Phi(a)$ for some $a\in A$, then $\phi_x(T)=\phi(ax)=\pi(a)\phi(x)=
\tilde{\pi}(T)\phi(x)=\psi_x(T)$ (since $\phi$ is a left $a$-module map), hence by 
weak* continuity
$\phi_x(T)=\psi_x(T)$ for all $T\in\tilde{A}$. This means that
$\tilde{\pi}(T)\phi(x)=\phi(\alpha(T)x)$, hence, since $E$ is a $C$-bimodule map fixing
elements of $C$,
$$E(\tilde{\pi}(T))\phi(x)=\phi(\alpha(T)x)\ \ \mbox{for all}\ T\in\tilde{A}\ \mbox{and}\ 
x\in\kh.$$ Using (\ref{34}) it follows now that
$$\phi(a_Tx)=\pi(a_T)\phi(x)=E(\tilde{\pi}(T))\phi(x)=\phi(\alpha(T)x),$$
hence $\phi((a_T-\alpha(T))\kh)=0$ and $a_T=\alpha(T)$ since the kernel of $\phi$ does
not contain any nontrivial right ideal.  This shows that $\alpha(\tilde{A})
\subseteq A$. But, since $\alpha$ is weak* continuous, it is a well known consequence of 
the Kaplansky density theorem (together with Alaoglu's theorem) that 
$\alpha(\tilde{A})=\overline{A}$.
Thus, $A=\overline{A}$.
\end{proof}

\begin{example}\label{e33} We note that the assumption in Theorem \ref{th32} that $C$ is
injective is not redundant. To see this, let $\h=\ell^2$, $A=c_0$ (the sequences, converging
to $0$ identified as diagonal operators relative to the standard orthonormal basis of $\ell^2$), $C=\kh$, $\pi:A\to C$ the inclusion
and $\phi:\kh\to C$ the identity mapping. Then the relative bicommutant of $A$ in $C$ is 
$A$, but $A$ is not a von Neumann algebra. 
\end{example}

Clearly, requiring merely that a C$^*$-algebra $A$ is equal to its relative bicommutant 
in an injective C$^*$-algebra containing $A$, does not imply that $A$ is a von Neumann
algebra (since any C$^*$-algebra $A$ is equal to its relative bicommutant in
$A$). But presently the author does not know if each monotone complete C$^*$-algebra
$A$ satisfies this bicommutation condition.

For C$^*$-algebras $A,B$ admitting faithful cyclic Hilbert modules, one could now deduce
from Theorems \ref{th32}, \ref{th31} and \ref{th21} that the pair $(A,B)$ admits a duality
(if and) only if $A$ and $B$ are von Neumann algebras. But we shall prove this conclusion
without assuming the existence of faithful cyclic Hilbert modules.

\begin{lemma}\label{le40}\cite{BZ} For any operator space $X$ and index set $\bj$ the multiplier algebra $\al(\matj{X})$
can naturally be identified (completely isometrically) as a subspace in $\matj{\al(X)}$.
\end{lemma}

\begin{proof}
First observe, from the fact that each left multiplier $\theta$ on $\matj{X}$ 
is a module map over $\matj{\bc}
\subseteq\ar(\matj{X})$, that $\theta$ must act on a 
matrix $[x_j]\in\matj{X}$, decomposed into columns 
$x_j\in\colj{X}$, as $\theta([x_j])=[\phi(x_j)]$, where $\phi\in{\rm CB}(\colj{X})$.
Moreover, using Theorem \ref{BEZ}(i) it follows that $\phi$ must be a left multiplier on
$\colj{X}$. Thus, a left multiplier on $\matj{X}$ is just a left multiplier on $\colj{X}$
applied to all columns of a matrix. The inclusion $\al(\colj{X})\subseteq\matj{\al(X)}$ 
is proved in
\cite[Remark following 5.10.1]{BZ}. To explain a slightly different 
approach, let $Y=I(X)$ be the injective envelope of $X$. Then $\colj{Y}$ is the injective
envelope of $\colj{X}$ by \cite[4.6.12]{BLM}. Assuming that $[YY^*]$ is contained in a
von Neumann algebra $\mathcal{R}$ (which may be taken of the form $\mathcal{R}=\bh$), we have
that the C$^*$-algebra $[\colj{Y}\rowj{Y^*}]$ is contained in $\matj{\mathcal{R}}$, 
hence the same holds for its multiplier C$^*$-algebra $\al(\colj{Y})$. But it follows from Theorem
\ref{BEZ}(ii) that $\al(\colj{X})\subseteq\al(\colj{Y})$ since $\colj{Y}$ is the injective
envelope of $\colj{X}$, hence $\al(\colj{X})\subseteq\matj{\mathcal{R}}$. Finally, a matrix
$[T_{ij}]\in\matj{\mathcal{R}}$ multiplies $\colj{X}$ into itself only if $T_{ij}X\subseteq X$
for all $i,j\in\bj$, hence $\al(\colj{X})\subseteq\matj{\al(X)}$. 
\end{proof}

\begin{lemma}\label{le42} Let  $A$ be  a 
C$^*$-subalgebra of $\al(X)$. Then for any index set $\bj$ the bicommutant 
of $A^{(\bj)}$ in $\al(\matj{X})$ satisfies
$$(A^{(\bj)})^{cc}=(A^{cc})^{(\bj)},$$
where $A^{cc}$ is the bicommutant of $A$ in $\al(X)$.
\end{lemma}

\begin{proof} By Lemma \ref{le40}
$\al(\matj{X})\subseteq\matj{\al(X)}$.
Moreover, $\al(\matj{X})$ contains the algebra $\mathcal{F}$ of all
finitely supported matrices in $\matj{\al(X)}$.
It follows that for a C$^*$-subalgebra $A$ of $\al(X)$, the commutant $(A^{(\bj)})^c$
of $A^{(\bj)}$ in $\al(\matj{X})$ contains the set
$$C:=\{\phi=[\phi_{ij}]\in\mathcal{F}:\ \phi a^{(\bj)}=a^{(\bj)}\phi\ \forall a\in A\}
=\{[\phi_{ij}]\in\mathcal{F}:\ \phi_{ij}\in A^c\}.$$
This implies that $(A^{(\bj)})^{cc}\subseteq C^c$. But a standard computation shows that
$C^c$ consists of diagonal matrices $x^{(\bj)}$, where $x\in A^{cc}$, hence
$(A^{(\bj)})^{cc}\subseteq(A^{cc})^{(\bj)}$. For the reverse inclusion, note that,
given $x\in A^{cc}$, $x^{(\bj)}$ commutes with all matrices $[\phi_{ij}]\in
(A^{(\bj)})^c=\matj{A^c}\cap\al(\matj{X})$.
\end{proof}

\begin{remark}\label{re61} If $Z$ is a cogenerator of an ample category $\om$ of
operator $A,B$-bimodules, then the $*$-homomorphism $A\to\al(Z)$, by which the
left $A$-module structure is introduced to $Z$, is injective. This can be deduced from 
Definition \ref{d41}(i) and Theorem \ref{th4}.
\end{remark}

\begin{theorem}\label{th43} Let $\om$ be an ample subcategory of $\aob$ and $Z$ an
injective cogenerator in $\om$. Suppose that $A$, regarded as a C$^*$-subalgebra of $\al(Z)$,
is equal to its bicommutant $A^{cc}$ in $\al(Z)$. Then $A$ is a W$^*$-algebra. 
\end{theorem}

\begin{proof} Let $\h\in\ha$ and $\k\in\hb$ be faithful Hilbert modules. By Theorem \ref{th4} 
there is a completely isometric $A,B$-bimodule
map $\phi:\bkh\to\matk{Z}$ for some index set $\bk$. By Theorem \ref{th31} $\phi$ induces
an  $A$-bimodule map $\pl:\bh\to\al(\matk{Z})$ such that the kernel of $\pl$ does not contain
any nontrivial right ideals. From the hypothesis and Lemma \ref{le42} the image
of $A$ in $\al(\matk{Z})$ is equal to its relative bicommutant, hence by Theorem \ref{th32}
$A$ is a von Neumann algebra on $\h$. 
\end{proof}

By Theorem \ref{th43} a pair of C$^*$-algebras $(A,B)$ admits a duality (in the sense defined
in the Introduction) only if $A$ and $B$ are von Neumann algebras. Conversely, each pair 
$(A,B)$ of von Neumann algebras admits a duality by the von Neumann
bicommutation theorem, since $\bkh$ is an injective cogenerator for normal operator 
$A,B$-bimodules if $\h$
and $\k$ are the universal normal Hilbert modules over $A$ and $B$, respectively. 

\section{Initial cogenerators}

In this Section $M$ and $N$ are von Neumann algebras, $\nhm$ the category of normal Hilbert
$M$-modules and $\mnon$ the category of all {\em normal}
operator $M,N$-bimodules. (We do not require that a normal bimodule $X$ is a dual space,
only that there exists a completely
isometric $M,N$-bimodule map from $X$ into $\bkh$ for some normal Hilbert modules $\h\in\nhm$ 
and $\k\in\nhn$.)  Recall from the Introduction that a cogenerator $Z$ in $\mnon$ is
{\em countably initial} if for each cogenerator $X$ there is a completely isometric
$M,N$-bimodule map from $Z$ into $X^{\bn}$. 

If $M$ and $N$ admit  cyclic  modules 
$\h\in\nhm$ and $\k\in\nhn$ such that all normal states on
$M$ and $N$ are vector states coming from vectors in $\h$ and $\k$, then $\bkh$ is
a  cogenerator in $\mnon$ (see the paragraph following Definition \ref{dc}). 
Moreover, it follows
from  Theorem \ref{th21} that $\bkh$ is countably initial. 
We shall show that the condition on $M$ and $N$  
is satisfied if (and only if)
the centers of  $M$ and $N$ 
are $\sigma$-finite. This will imply the first part of the following Theorem. 

\begin{theorem}\label{th61} If the centers of $M$ and $N$ are $\sigma$-finite,
then $\mnon$ has a countably initial cogenerator
of the form $\bkh$, where $\h\in\nhm$ and $\k\in\nhn$ are cyclic and
such that all normal states on $M$ and $N$ are vector states coming
from vectors in $\h$ and $\k$, respectively. If the center of $M$ 
is not $\sigma$-finite, then ${{_{M}{\rm NOM}_{\bc}}}$ has no countably initial cogenerators.
\end{theorem}

By the observation above, one direction of the Theorem is an immediate consequence of 
Theorem \ref{th21} and the following lemma.

\begin{lemma}\label{le} If the center of  $M$ is $\sigma$-finite, then there
exists a  cyclic  $\h\in\nhm$ such that all normal
states on $M$ are vector states coming from vectors in $\h$.
\end{lemma}

Since we have not found a reference for the lemma, we include a proof.
\begin{proof}[Proof of Lemma \ref{le}] We shall use repeatedly the fact that a sum of a countable family 
of centrally orthogonal cyclic projections in a von Neumann algebra is cyclic \cite[5.5.10]{KR}.
Decomposing $M$
into the direct sum of the finite, the properly infinite (and semifinite) and the purely infinite
part, we may consider each part separately. 

If $M$ is finite, then $M$ is $\sigma$-finite
(since the center of $M$ is) by \cite[8.2.9]{KR} or \cite[V.2.9]{T}, hence it can be represented on a Hilbert 
space $\h$ so that
it has a cyclic and separating (trace) vector and then all normal states on $M$
are vector states coming from vectors on $\h$ by \cite[7.2.3]{KR} (or \cite[V.1.12]{T}). 

If $M$ is properly infinite and semifinite,
then by \cite[V.1.40]{T} $M$ is a direct sum of a countable family of algebras of the form 
$N\overline{\otimes}\B(\g)$, where $N$ is finite and $\g$ is infinite-dimensional. Representing
$N$ on a space $\h_N$ where it has a cyclic and separating vector and taking the tensor 
product of this representation with  the countable multiple of the identity representation of 
$\B(\g)$ (cyclic), gives a 
cyclic representation
of $N\overline{\otimes}\B(\g)$ on $\h_N\otimes\g\otimes\ell^2_{\bn}$ such that all normal states are vector states. 

Finally, if
$M$ is purely infinite, acting on a Hilbert space $\g$, let 
$\{\com{e}_j:\ j\in\bj\}$ be a maximal family of cyclic projections in $\com{M}$
with mutually orthogonal central carriers $p_j$. Then $\bj$ is countable and,
since the central carrier of the projection $\com{e}:=\sum_{j\in\bj}\com{e}_j$ is $1$
by maximality,
the representation of $M$ on $e^{\prime}\g$ ($a\mapsto a\com{e}$) is faithful
and cyclic. In particular $\com{(M\com{e})}$ has a separating vector and is therefore
$\sigma$-finite. However, since $M$ is of type III with $\sigma$-finite center,
by \cite[V.3.2]{T} $M$ has a unique faithful normal
representation with $\sigma$-finite commutant, up to a unitary equivalence. So, 
the countable multiple of the representation
of $M$ on $\com{e}\g$ is unitarily equivalent to the same representation and therefore
each normal state on $M$ must be a vector state coming from a vector in $\h:=\com{e}\g$. 
\end{proof}

We remark that the Hilbert modules $\h$ and $\k$ in Theorem \ref{th61} are not
necessarily the ones on which $M$ and $N$ are in the standard form.
To see this, we may consider $M=N=\B(\ell^{2}_{\bi})=\mati{\bc}$, where $\bi$ is an uncountable set.
Then $X:=\B(\ell^2_{\bi}\otimes\ell^2_{\bi})=\mati{\B(\ell^2_{\bi})}$ is an injective
cogenerator for operator $M$-bimodules, but can not be embedded into $\B(\ell^2_{\bi}
\otimes\ell^2_{\bn})^{\bn}=\mati{\B(\ell^2_{\bn})}^{\bn}=\mati{\B(\ell^2_{\bn})^{\bn}}$ 
as an $M$-bimodule. Namely, such an
embedding would  be the amplification of an embedding of $\B(\ell^2_{\bi})$ into
$\B(\ell^2_{\bn})^{\bn}$, which can not exist if $\bi$ is large enough.

Since a representation (state) of a C$^*$-algebra $A$ is just a normal representation
(state) of its universal von Neumann envelope $\tilde{A}$ and the center of $\tilde{A}$
is $\sigma$-finite if and only if each family of disjoint cyclic representations of $A$
is countable (this follows using \cite[10.3.3]{KR}), the method of Theorem \ref{th61} 
proves also the first part of the following Proposition.

\begin{proposition}\label{pr0} If all families of disjoint cyclic representations of 
C$^*$-algebras $A$ and $B$ are countable, then the category $\aob$ has an initial
strict cogenerator of the form $\bkh$, where $\h\in\ha$ and $\k\in\hb$ are cyclic.
If $A$ has an uncountable family of disjoint (cyclic) representations, then
${{_{A}{\rm OM}_{\bc}}}$ has no initial strict cogenerators.
\end{proposition}

It is known that separable non type I C$^*$-algebras have uncountably many disjoint
(cyclic) representations \cite[6.8.5]{Pe}. Thus, separable C$^*$-algebras satisfying the condition of the above
Proposition, are of type I and have only countably many inequivalent irreducible representations.

We still have to prove the last parts of Theorem \ref{th61} and Proposition \ref{pr0}.
Let $\{p_i:\ i\in\bi\}$ be a (fixed through the rest of the section) maximal orthogonal family of central projections in $M$ that
are $\sigma$-finite in the center. Hence $\sum_{i\in\bi}p_i=1$ by maximality. We shall assume
that $\bi$ is uncountable (otherwise the center of $M$ is $\sigma$-finite) and denote by 
$\bs$ the family of all countable subsets of $\bi$.
For each $s\in\bs$ let $p_s:=\sum_{i\in s}p_i$.
By Lemma \ref{le} we may choose for each $i\in\bi$ a normal cyclic Hilbert $p_iM$-module $\h_i$ such that
all normal states on $p_iM$ are vector states from vectors in $\h_i$. Let $\h_M=\oplus_i\h_i$, so that
$\h_i=p_i\h_M$. Let $\k$ be any fixed Hilbert space ($\k=\bc$ will be sufficient for our
application here).  

\begin{definition}\label{d6}
Let $Z_M$ be the $M$-submodule
of $\B(\k^{\bs},\h_M)={\rm R}_{\bs}(\B(\k,\h_M)$, consisting of all $z=[z_s]_{s\in\bs}$,
where $z_s\in\B(\k,\h_M)$ are such that:

(i) $z_s=p_sz_s$ for all $s\in\bs$ and

(ii) for each $i\in\bi$ the set of all $s\in\bs$ such that $p_iz_s\ne0$ is countable.

\end{definition}

Each normal cyclic Hilbert $M$-module $\h$ is the orthogonal sum of submodules $p_i\h$,
where the set $s$ of all $i\in\bi$ such that $p_i\h\ne0$ is  countable (by cyclicity). 
Since each $p_i\h$ is cyclic over $p_iM$ and all normal states on $p_iM$ are vector states
coming from vectors in $p_i\h_M$, $p_i\h$ is
contained in $p_i\h_M$ (up to a unitary equivalence), hence $\h$ is contained in 
$\oplus_{i\in s}p_i\h_M=p_s\h_M$. Since $p_s\h_M$ is contained in $\h_M$ and in $Z_M$
(if $\k=\bc$ then $p_s\h_M$ is just the $s$-th column of $Z_M$), it follows now from
Proposition \ref{pr41} that $\h_M$ and $Z_M$ are cogenerators in ${{_M{\rm NOM}_{\bc}}}$.

\begin{lemma}\label{le61} For each $j\in\bi$ let $y_j=p_jy_j$ be a norm $1$ element in
$Z_M^{\bn}$ and let $i\in\bi$ be fixed. Then  the set $\bj$ of all 
$j\in\bi$ such
that
$$ \|\left[\begin{array}{c}
y_i\\
y_j
\end{array}\right]\|>1
$$
is countable.
\end{lemma}

\begin{proof} Let $y_{j,n}\in Z_M$ be the components of $y_j$. If $\bj$ is 
uncountable, then for some $n$ the set $\bj_0$ of all $j\in\bj$ such that
\begin{equation}\label{e7}\|\left[\begin{array}{c}
y_{i,n}\\
y_{j,n}
\end{array}\right]\|>1
\end{equation}
is uncountable.
Put $z(j)=y_{j,n}$  and let $z(j)_s$ ($\in\B(\k,\h_M)$) be the components of 
$z(j)=p_jz(j)\in{\rm R}_{\bs}(\B(\k,\h_M))$. 
Let $\bs_i=\{s\in\bs:\ z(i)_s\ne0\}$ (a countable set by Definition \ref{d6}). 

For each $z=[z_s]\in Z_M$ (where $z_s\in\B(\k,\h_M)$ are the components
of $z$) we have by definition of $Z_M$ that $z_s=p_sz_s$, hence $p_jz_s=0$ if $j\notin s$.
Thus, if $j$ is outside of the countable set $\bi_0:=
\cup_{s\in\bs_i}s$, then  $p_jz_s=0$ for all $s\in\bs_i$. This holds in particular for 
$z=z(j)$, hence, if $j\notin\bi_0$, then $p_jz(j)_s=0$ for all $s\in\bs_i$.
Since $z(j)=p_jz(j)$, we also have that $z(j)_s=p_jz(j)_s$ for all $s\in\bs$,
and it follows that $z(j)_s=0$ if $s\in\bs_i$ and $j\notin\bi_0$. Thus, if $j\notin\bi_0$, then the
two rows $z(i)=[z(i)_s]_{s\in\bs}$ and $z(j)=[z(j)_s]_{s\in\bs}$ have disjoint supports and
therefore 
$$\norm{\left[\begin{array}{c}
z(i)\\
z(j)
\end{array}\right]}=\max\{\norm{z(i)},\norm{z(j)}\}=1.$$
It follows that the inequality (\ref{e7}) can hold only if $j\in\bi_0$. In other words,
$\bj_0\subseteq\bi_0$, hence $\bj_0$ is countable, a contradiction.
\end{proof}

\begin{proof}[Proof of Theorem \ref{th61}] It only remains to prove that the category
${{_M{\rm NOM}_{\bc}}}$ has no countably initial cogenerators. 
Suppose the contrary, that $Z_0$ is such a cogenerator. Choose a small $\varepsilon>0$
($\varepsilon<\sqrt{2}-1$). We shall use the notation introduced in
and above Definition \ref{d6}. Since each $p_j\h_M$ is cyclic,
by Proposition \ref{pr41} $Z_0^{\bn}$ contains $p_j\h_M$ completely isometrically as an $M$-module
for each $j\in\bi$, hence  we can choose an element
$y_j=p_jy_j\in Z_0^{\bn}$ of norm $1$. Since $\h_M$ is a cogenerator and $Z_0$ is countably 
initial, we may regard $Z_0$ as an $M$-submodule in $\h_M^{\bn}$, hence $Z_0^{\bn}$ as
a submodule in $\h_M^{\bn\times\bn}\cong \h_M^{\bn}$. Let $y_{j,n}\in\h_M$ be the components
of $y_j$ regarded as an element in $\h_M^{\bn}$. Then for each $j\in\bi$ there exists 
$n\in\bn$ such that $\norm{y_{j,n}}>1-\varepsilon/2$.
Since $\bi$ is uncountable, there exists a $n$ such that $\norm{y_{j,n}}>1-\varepsilon/2$
for all $j$ in an uncountable subset $\bj$ of $\bi$.  
If $i,j\in\bj$, then, since  $y_{i,n}$ and $y_{j,n}$ are  vectors in the Hilbert space $\h_M$, 
this implies that 
$$\|\left[\begin{array}{c}
y_{i,n}\\
y_{j,n}
\end{array}\right]\|>\sqrt{2}-\varepsilon\ \ \mbox{if}\ \ i,j\in\bj.$$
It follows that 
\begin{equation}\label{65}\|\left[\begin{array}{c}
y_i\\
y_j
\end{array}\right]\|>\sqrt{2}-\varepsilon\ \ \mbox{if}\ \ i,j\in\bj.
\end{equation}

Since $Z_0$ is countably initial and $Z_M$ is a cogenerator, $Z_0$ is contained as
an $M$-submodule in $Z_M^{\bn}$ completely isometrically, hence $Z_0^{\bn}$ is 
contained in $Z_M^{\bn\times\bn}\cong Z_M^{\bn}$. But then by Lemma \ref{le61}
each set of elements $y_j$ in $Z_0^{\bn}$ satisfying (\ref{65}) (for a fixed $i$)  must be countable,
which contradicts the fact that $\bj$ in (\ref{65}) is uncountable.
\end{proof}

\begin{proof}[Proof of Proposition \ref{pr0}] We still have to prove that 
${{_A{\rm OM}_{\bc}}}$
has no initial strict cogenerators if $A$ has uncountably many disjoint cyclic representations.
The condition on $A$ means that the center of the universal von Neumann envelope $\tilde{A}$
of $A$ is not $\sigma$-finite. If we put $M=\tilde{A}$ and construct $\h_M$ and $Z_M$ as above,
then $\h_M$ and $Z_M$ are strict cogenerators in ${{_A{\rm OM}_{\bc}}}$ by Proposition
\ref{pr41}. Thus, if there exists an initial strict cogenerator $Z_0$, then both, $\h_M$ and $Z_M$
contain a copy of $Z_0$. Similarly as in the above proof of Theorem \ref{th61} (but easier), 
the fact that $Z_0$ is a strict cogenerator implies the existence of an uncountable set 
$\{y_j:\ j\in\bi\}$ satisfying (\ref{65}), while the inclusion $Z_0\subseteq Z_M$ together
with Lemma \ref{le61} shows that this is impossible.
\end{proof}

The same technique shows that if the center of at least one of the algebras $M$ or $N$ is not 
$\sigma$-finite
the category $\mnon$ has no countably initial cogenerators of the form $\bkh$. It seems
natural to conjecture that in this case $\mnon$ has no countably initial cogenerators at all.
The above technique can be upgraded to prove this in the case one
of the two algebras  has a separable predual (while the center of the other is not $\sigma$-finite),
but the general case remains open.

\end{document}